\documentclass[12pt,draft]{amsart}
\usepackage{amsmath,amsthm,latexsym,amscd,amsbsy,amssymb,fleqn,leqno}
\chardef\bslash=`\\ 

\makeatletter
\def\verbatim{\interlinepenalty\@M \@verbatim
  \leftskip\@totalleftmargin\advance\leftskip2pc
  \frenchspacing\@vobeyspaces \@xverbatim}
\makeatother
\makeatletter
  \def\dgt@k{\dg@DX=-3 \dg@DY=2 \dg@SIZE=3} 
\makeatother

\makeatletter
  \def\dgt@kk{\dg@DX=3 \dg@DY=-1 \dg@SIZE=3}%
\makeatother

\theoremstyle{plain}
\newtheorem{thm}{Theorem}[section]
\newtheorem{cor}[thm]{Corollary}
\newtheorem{lem}[thm]{Lemma}

\theoremstyle{definition}

\numberwithin{equation}{section}

\newcounter{rmnum}

\def\symbolnote#1#2{\let\thefootn=\thefootnote%
\renewcommand{\thefootnote}{\fnsymbol{footnote}}%
\footnotemark[#1]%
\footnotetext[#1]{#2}%
\let\thefootnote=\thefootn
}

\newfont{\bbb}{msbm10 scaled \magstep1}
\newfont{\bbc}{msbm8 scaled \magstep0}

\newcommand{\R}{\mbox{\bbb R}}

\newcommand{\N}{\mbox{\bbb N}}

\newcommand{\uin}{\mbox{\bbb I}}

\begin{document}


\title[On finite-dimensional  maps]{On finite-dimensional maps II}

\author{H. Murat Tuncali}
\address{Department of Mathematics,
Nipissing University,
100 College Drive, P.O. Box 5002, North Bay, ON, P1B 8L7, Canada}
\email{muratt@nipissingu.ca}
\thanks{The first author was partially supported by NSERC grant.}

\author{Vesko Valov}
\address{Department of Mathematics, Nipissing University,
100 College Drive, P.O. Box 5002, North Bay, ON, P1B 8L7, Canada}
\email{veskov@nipissingu.ca}
\thanks{The second author was partially supported by Nipissing University Research Council Grant.}

 

\begin{abstract}{Let $f\colon X\to Y$ be a perfect $n$-dimensional surjective map of paracompact spaces and $Y$  a $C$-space. We consider the following property of continuous maps  $g\colon X\to\uin^k=[0,1]^k$, where $1\leq k\leq\omega$:  each $g(f^{-1}(y))$, $y\in Y$, is at most $n$-dimensional. It is shown that all maps $g\in C(X,\uin^{n+1})$ with the above property  form a dense $G_{\delta}$-set in the function space
$C(X,\uin^{n+1})$  equipped with the source limitation topology. Moreover, for every $n+1\leq m\leq\omega$ the space $C(X,\uin^m)$ contains a dense $G_{\delta}$-set of maps having this property.}
\end{abstract}

\maketitle
\markboth{H. M.~Tuncali and V.~Valov}{On finite-dimensional maps}


\section{Introduction}

This note is inspired by a result of Uspenskij \cite[Theorem 1]{vu:01}. Answering a question of R. Pol, Uspenskij proved the following theorem: Let $f\colon X\to Y$ be a light map (i.e. every fiber $f^{-1}(y)$ is $0$-dimensional) between compact spaces and $\mathcal{A}$ be the set of all functions $g\colon X\to\uin=[0,1]$ such that $g(f^{-1}(y))$ is 0-dimensional for all $y\in Y$. Then $\mathcal{A}$ is a dense $G_{\delta}$-subset of the function space $C(X,\uin)$ provided $Y$ is a $C$-space (the case when $Y$ is countable-dimensional was established earlier by Torunczyk). We extend this result as follows:  

\begin{thm}
Let $f\colon X\to Y$ be a $\sigma$-perfect surjection such that $\dim f\leq n$ and $Y$ is a paracompact $C$-space.  Let 
$\mathcal{H}=\{g\in C(X,\uin^{n+1}): \dim g(f^{-1}(y))\leq n\mbox{ for each $y\in Y$}\}$. Then $\mathcal{H}$ is dense and 
$G_{\delta}$ in $C(X,\uin^{n+1})$ with respect to the source limitation topology.
\end{thm}

\begin{cor}
Let $X$, $Y$ and $f$ satisfy the hypotheses of Theorem $1.1$ and $n+1\leq m\leq\omega$. Then, there exists a dense $G_{\delta}$-subset $\mathcal{H}_m$ of $C(X,\uin^m)$ with respect to the source limitation topology such that $\dim g(f^{-1}(y))\leq n$ for every $g\in\mathcal{H}_m$ and $y\in Y$.
\end{cor}

Here, $\dim f=\sup\{\dim f^{-1}(y):y\in Y\}$ and $f$ is said to be $\sigma$-perfect if there exists a sequence $\{X_i\}$ of closed subsets of $X$ such that each restriction map $f|X_i$ is perfect and the sets $f(X_i)$ are closed in $Y$. 
The $C$-space property was introduced by Haver \cite{wh:74} for compact metric spaces and then extended by Addis and Gresham \cite{ag:78} for general spaces (see \cite{re:95} for the definition and some properties of $C$-spaces). Every countable-dimensional (in, particular, every finite-dimensional) paracompact space has property $C$, but there exists a compact metric $C$-space which is not countable-dimensional
\cite{rp:81}. For any
spaces $X$ and $Y$ by $C(X,Y)$ we denote the set of all continuous maps from
$X$ into $Y$. If $(Y,d)$ is a metric space, then the source limitation topology on $C(X,Y)$ is defined in the following way:
a subset $U\subset C(X,Y)$ is open in $C(X,Y)$ with respect to
the source limitation topology provided for every $g\in U$ there exists 
a continuous function $\alpha\colon X\to (0,\infty)$ such that $\overline{B}(g,\alpha)\subset U$, where $\overline{B}(g,\alpha)$ denotes the set  
$\{h\in C(X,Y):d(g(x),h(x))\leq\alpha (x)\hbox{}~~\mbox{for each 
$x\in X$}\}$. 
The source limitation topology is also known as the fine topology and $C(X,Y)$ with this topology has Baire property provided $(Y,d)$ is a complete metric space \cite{jm:75}.  Moreover, the source limitation topology on $C(X,Y)$ doesn't depend on the metric of $Y$ when $X$ is paracompact \cite{nk:69}. 

All single-valued maps under discussion are continuous, and all function spaces, if not explicitely stated otherwise, are equipped with the source limitation topology.

\section{Proofs}

Let show first that the proof of Theorem 1.1 can be reduced to the case when $f$ is perfect.  Indeed, we fix a sequence $\{X_i\}$ of closed subsets of $X$ such that each map $f_i=f|X_i\colon X_i\to Y_i=f(X_i)$ is perfect and $Y_i\subset Y$ is closed. Consider the maps $\pi_i\colon C(X,\uin^{n+1})\to C(X_i,\uin^{n+1})$ defined by $\pi (g)=g|X_i$ and the sets 
$\mathcal{H}_i=\{g\in C(X_i,\uin^{n+1}): \dim g(f_i^{-1}(y))\leq n\mbox{ for each $y\in Y_i$}\}$.  If Theorem 1.1 holds for perfect maps, then every $\mathcal{H}_i$ is dense and $G_{\delta}$ in $C(X_i,\uin^{n+1})$, so are the sets $\pi_i^{-1}(\mathcal{H}_i)$ in $C(X,\uin^{n+1})$ because $\pi_i$ are open and surjective maps.  Finally, observe that $\mathcal{H}$ is the intersection of all $\mathcal{H}_i$ and since $C(X,\uin^{n+1})$ has Baire property, we are done.
 
Everywhere in this section $X$, $Y$, $f$ and $\mathcal{H}$ are fixed and satisfy the hypotheses of Theorem 1.1 with $f$ being perfect. Any finite-dimensional cube $\uin^k$ is considered with the Euclidean metric. We say that a set-valued map $\theta\colon H\to\mathcal{F}(Z)$, where $\mathcal{F}(Z)$ denotes the family of all closed subsets of the space $Z$, is upper semi-continuous (br. u.s.c.) if $\{y\in H:\theta(y)\subset W\}$ is open in $H$ for every open $W\subset Z$.  
In the above notation, $\theta$ is called lower semi-continuous if $\{y\in H:\theta(y)\cap W\neq\emptyset\}$ is open in $H$ whenever $W$ is open in $Z$.

\medskip
\noindent{\em Proof of Theorem $1.1$}.

For every open set $V$ in $\uin^{n+1}$ let $\mathcal{H}_{V}$ be the set of all $g\in C(X,\uin^{n+1})$ such that $V$ is not contained in any $g(f^{-1}(y))$, $y\in Y$.  Following the Uspenskij idea from \cite{vu:01}, it suffices to show that each set 
$\mathcal{H}_{V}$ is dense and open in $C(X,\uin^{n+1})$. Indeed, choose a countable base $\mathcal{B}$ in 
$\uin^{n+1}$. Since a subset of $\uin^{n+1}$ is at most $n$-dimensional if and only if it doesn't contain any $V\in\mathcal{B}$, we have that
$\mathcal{H}$ is the intersection of all $\mathcal{H}_V$, $V\in\mathcal{B}$. But $C(X,\uin^{n+1})$ has the Baire property, so
$\mathcal{H}$ is dense and $G_{\delta}$ in $C(X,\uin^{n+1})$ and we are done.

\begin{lem}
The set $\mathcal{H}_{V}$ is open in $C(X,\uin^{n+1})$ for every open $V\subset\uin^{n+1}$.
\end{lem}

\begin{proof}
Fix an open set $V$ in $\uin^{n+1}$ and $g_0\in\mathcal{H}_V$.   We are going to find a continuous function $\alpha\colon X\to (0,\infty)$ such that $\overline{B}(g_0,\alpha)\subset\mathcal{H}_V$. To this end, let $p\colon Z\to Y$ be a perfect surjection with $\dim Z=0$ and define $\psi\colon Y\to\mathcal{F}(\uin^{n+1})$ by
$\psi (y)=g_0(f^{-1}(y))$, $y\in Y$. Since $f$ is perfect,  $\psi$ is upper semi-continuous 
and compact-valued.
Now, consider the set-valued map $\psi_1\colon Z\to\mathcal{F}(\uin^{n+1})$, 
$\psi_1=\psi\circ p$. Obviously, $g_0\in\mathcal{H}_V$ implies $\overline{V}\backslash\psi_1(z)\neq\emptyset$ for every $z\in Z$. Moreover, $\psi_1$ is also upper semi-continuous, in particular it has a closed graph.  Then, by a result of Michael \cite[Theorem 5.3]{em:88}, there exists a continuous map $h\colon Z\to\uin^{n+1}$ such that $h(z)\in\overline{V}\backslash\psi_1(z)$, $z\in Z$. Next, 
consider the u.s.c. compact-valued map $\theta\colon Y\to\mathcal{F}(\uin^{n+1})$,
$\theta (y)=h(p^{-1}(y))$, $y\in Y$. We have
$\emptyset\neq\theta(y)\subset\overline{V}$ and $\theta(y)\cap\psi(y)=\emptyset$ for all $y\in Y$. Hence, the function $\alpha_1\colon Y\to\R$,
$\alpha_1(y)=d(\theta(y),\psi(y))$, is positive, where $d$ is the Euclidean metric on $\uin^{n+1}$. Since, both $\theta$ and $\psi$ are upper semi-continuous, $\alpha_1$ has the following property: $\alpha_1^{-1}(a,\infty)$ is open in $Y$ for every $a\in\R$.  Finally, take a continuous function $\alpha_2\colon Y\to (0,\infty)$ with $\alpha_2(y)<\alpha_1(y)$ for every $y\in Y$ (see, for example, \cite{jd:44}) and define $\alpha=\alpha_2\circ f$. It remains to observe that, if $g\in\overline{B}(g_0,\alpha)$ and $y\in Y$, then 
$\theta(y)\subset\overline{V}\backslash g(f^{-1}(y))$. So, $g(f^{-1}(y))$ doesn't contain $V$ for all $y\in Y$, i.e. $\overline{B}(g_0,\alpha)\subset\mathcal{H}_V$.
\end{proof}

{\bf Remark.} Analyzing the proof of Lemma 2.1, one can see that we proved the following more general statement:  Let $h\colon\overline{X}\to\overline{Y}$ be a perfect surjection between paracompact spaces and $K$ a complete metric space. Then, for every open $V\subset K$ the set of all maps $g\in C(\overline{X},K)$  with $V\not\subset g(h^{-1}(y))$ for any $y\in\overline{Y}$ is open in $C(\overline{X},K)$.

The remaining part of this section is devoted to the proof that each $\mathcal{H}_{V}$ is dense in $C(X,\uin^{n+1})$, which is finally accomplished by Lemma 2.6.  

\begin{lem}
Let $Z$ and $K$ be compact spaces and $K_0=\cup_{i=1}^{\infty}K_i$ with each $K_i$ being a closed $0$-dimensional subset of $K$. Then the set 
$\mathcal{A}=\{g\in C(Z\times K,\uin)$: $\dim g(\{z\}\times K_0)=0\mbox{ for every $z\in Z$}\}$ is dense and $G_{\delta}$ in $C(Z\times K,\uin)$. 
\end{lem}

\begin{proof}
Since, for every $i$, the restriction map $p_i\colon C(Z\times K,\uin)\to C(Z\times K_i,\uin)$ is a continuous open surjection, we can assume that $K_0=K$ and $\dim K=0$. Then $\mathcal{A}$ is the intersection of the sets $\mathcal{A}_V$, $V\in\mathcal{B}$, 
where  $\mathcal{B}$ is a countable base of $\uin$ and $\mathcal{A}_V$
consists of all
$g\in C(Z\times K,\uin)$ such that $V\not\subset g(\{z\}\times K)$ for every $z\in Z$.   By the remark after Lemma 1.1, every $\mathcal{A}_V\subset C(Z\times K,\uin)$ is open, so $\mathcal{A}$ is $G_{\delta}$. 
It remains only to show that $\mathcal{A}$ is dense in $C(Z\times K,\uin)$.  
Since $K$ is $0$-dimensional, the set $C_K=\{h\in C(K,\R): h(K)\mbox{ is finite}\}$ is dense in $C(K,\R)$. Hence, by the Stone-Weierstrass theorem, all polynomials of elements of the family $\gamma=\{t\cdot h:t\in C(Z,\R), h\in C_K\}$ form a  dense subset $\mathcal P$ of $C(Z\times K,\R)$. We fix a retraction $r\colon\R \to\uin$ and define 
$u_r\colon C(Z\times K,\R)\to C(Z\times K,\uin)$, $u_r(h)=r\circ h$. Then $u_r(\mathcal P)$ is dense in $C(Z\times K,\uin)$.
It is easily seen that  every $g\in u_r(\mathcal P)$ has the following property: $g(\{z\}\times K)$ is finite for every $z\in Z$. So, 
$u_r(\mathcal P)\subset\mathcal A$, i.e. $\mathcal{A}$ is dense in $C(Z\times K,\uin)$. 
\end{proof}

\begin{lem}
Let $M$ and $K$ be compact spaces with $\dim K\leq n$ and $M$ metrizable. If $V\subset\uin^{n+1}$ is open, then the set of all maps $g\in C(M\times K,\uin^{n+1})$ such that  $V\not\subset g(\{y\}\times K)$  for each $y\in M$ is dense in $C(M\times K,\uin^{n+1})$. 
\end{lem}

\begin{proof}
We are going to prove this lemma by induction with respect to the dimension of $K$. According to Lemma 2.2, it is true if $\dim K=0$.  Suppose the lemma holds for any $K$ with $\dim K\leq m-1$ for some $m\geq 1$ and let $K$ be a fixed compact space with $\dim K=m$.   For 
$g^0\in C(M\times K,\uin^{m+1})$ and $\epsilon>0$ we need to find a function $g\in C(M\times K,\uin^{m+1})$ which is $\epsilon$-close to $g^0$ and $V\not\subset g(\{y\}\times K)$ for every $y\in M$. If $K$ is not metrizable, we represent it 
as the limit space of a $\sigma$-complete inverse system $\mathcal S=\{K_{\lambda},p_{\lambda}^{\lambda +1}:\lambda\in\Lambda\}$ such that each $K_{\lambda}$ is a metrizable compactum with $\dim K_{\lambda}\leq m$. Then $M\times K$ is the limit of the system $\{M\times K_{\lambda},id\times p_{\lambda}^{\lambda +1}:\lambda\in\Lambda\}$, where $id$ is the identity map on $M$.
Applying standard inverse spectra arguments (see \cite{book}), we can find $\lambda (0)\in\Lambda$ and $g_{\lambda (0)}\in C(M\times K_{\lambda(0)},\uin^{m+1})$ such that $g_{\lambda(0)}\circ (id\times p_{\lambda(0)})=g^0$, where $p_{\lambda(0)}\colon K\to K_{\lambda(0)}$ denotes the $\lambda(0)$-th limit projection of $\mathcal{S}$. 
Therefore, the proof is reduced to the case when $K$ is metrizable. 

Let $K$ be metrizable and $K=K_1\cup K_2$ such that $K_1$ is a $0$-dimensional $\sigma$-compact subset of $K$ and $\dim K_2\leq m-1$ (this is possible because $K$ is metrizable and $m$-dimensional, see \cite{re:95}). Let $g^0=g_1^0\times g_2^0$, where every $g_1^0$ is a function from $K$ into $\uin$, and $g_2^0\colon X\to\uin^m$.  We can assume that $V=V_1\times V_2$ with both $V_1\subset\uin$ and $V_2\subset\uin^m$ open.
According to Lemma 2.2, there exists a function $g_1\colon M\times K\to\uin$ which is $\displaystyle\frac{\epsilon}{\sqrt{2}}$-close to $g_1^0$ and such that $\dim g_1(\{y\}\times K_1)=0$ for every $y\in M$. Hence, 
$V_1$ is not contained in any of the sets $g_1(\{y\}\times K_1)$, $y\in M$. 

\medskip
{\em Claim. There exists an open set $A_1\subset K$ containing $K_1$ such that  
$\overline{V_1}\not\subset g_1(\{y\}\times A_1)$ for any $y\in M$}. 

\medskip
To prove the claim, we represent $K_1$ as the union of countably many compact $0$-dimensional sets $K_{1i}$ and consider the upper semi-continuous compact-valued maps $\psi_i\colon M\to\mathcal{F}({\uin})$ defined by $\psi_i(y)=g_1(\{y\}\times K_{1i})$. As in the proof of Lemma 2.1, we fix a $0$-dimensional space $Z$, a surjective perfect map $p\colon Z\to M$ and define the set-valued maps $\overline{\psi}_i\colon Z\to\mathcal{F}(\uin)$,  $\overline{\psi}_i=\psi_i\circ p$. It follows from our construction that each 
$\overline{\psi}_i(z)$, $z\in Z$, $i\in\N$, is $0$-dimensional. By \cite[Theorem 5.5]{em:88} (see also \cite[Theorem 1.1]{gv:01}), there is $h\in C(Z,\uin)$ such that 
$h(z)\in\overline{V_1}\backslash\cup_{i=1}^{\infty}\overline{\psi}_i(z)$, $z\in Z$. Then
$\theta\colon M\to\mathcal{F}(\uin)$, 
$\theta(y)=h(p^{-1}(y))$, is u.s.c. with 
$\emptyset\neq\theta(y)\subset\overline{V_1}\backslash g_1(\{y\}\times K_1)$ for every $y\in M$. Since the graph $G_{\theta}$ of $\theta$ is closed in $M\times\uin$, the set 
$U=\{(y,x)\in M\times K:(y,g_1(y,x))\not\in G_{\theta}\}$ is open in $M\times K$ and contains 
$M\times K_1$. So, $A_1=\{x\in K:M\times\{x\}\subset U\}$ is open in $K$ and contains $K_1$. Moreover, $\theta(y)\subset\overline{V_1}\backslash g_1(\{y\}\times A_1)$ for every $y\in M$,
which completes the proof of the claim.

\smallskip
Now, let $A_2=K\backslash A_1$. Obviously, $A_2$ is a compact subset of $K_2$, so $\dim A_2\leq m-1$.  According to the assumption that the lemma is true for any space of dimension $\leq m-1$, there exists
a map
$h_2\colon M\times A_2\to\uin^{m}$ which is $\displaystyle\frac{\epsilon}{\sqrt{2}}$-close to $g_2^0|(M\times A_2)$ and such that 
$\overline{V_2}\not\subset h_2(\{y\}\times A_2)$ for any $y\in M$.  We finally extend $h_2$ to a map $g_2\colon M\times K$ such that $g_2$ is $\displaystyle\frac{\epsilon}{\sqrt{2}}$-close to $g_2^0$.  Hence,

\medskip\noindent   
(1) \hbox{}~~~~~~$K=A_1\cup A_2$ and 
$\overline{V_j}\not\subset g_j(\{y\}\times A_j)$ for any $y\in M$, $j=1,2$. 

\medskip\noindent
Then the map $g=g_1\times g_2\colon M\times K\to\uin^{m+1}$ is $\epsilon$-close to $g_0$. It follows from (1) that $\overline{V}\not\subset g(\{y\}\times K)$ for any $y\in M$.  
 
\end{proof}

For any open $V\subset\uin^{n+1}$ we consider the set-valued map $\psi_V$ from $Y$ into $C(X,\uin^{n+1})$, given by 
$\psi_V(y)=\{g\in C(X,\uin^{n+1}): V\subset g(f^{-1}(y))\}$, $y\in Y$.

\begin{lem}
If $V\subset\uin^{n+1}$ is open and
$C(X,\uin^{n+1})$ is equipped with the uniform convergence topology, then   $\psi_V$ has a closed graph.  
\end{lem}

\begin{proof} 
Let $G_V\subset Y\times C(X,\uin^{n+1})$ be the graph of $\psi_V$ and $(y_0,g_0)\not\in G_V$.  Then $g_0\not\in\psi_V(y_0)$, so $g_0(f^{-1}(y_0))$ doesn't contain $V$. Consequently, there exists $z_0\in V\backslash g_0(f^{-1}(y_0))$ and let $\epsilon =d(z_0,g_0(f^{-1}(y_0)))$. Since $f$ is a closed map, there exists a neighborhood $U$ of $y_0$ in $Y$ with $d(z_0,g_0(f^{-1}(y)))>2^{-1}\epsilon$ for every $y\in U$. It is easily seen that $U\times B_{4^{-1}\epsilon}(g_0)$ is a neighborhood of $(y_0,g_0)$ in $Y\times C(X,\uin^{n+1})$ which doesn't meet $G_V$ (here $B_{4^{-1}\epsilon}(g_0)$ is the $4^{-1}\epsilon$-neighborhood of $g_0$ in $C(X,\uin^{n+1})$ with the uniform metric).
 Therefore $G_V\subset Y\times C(X,\uin^{n+1})$ is closed. 
\end{proof}

Recall that a closed subset $F$ of the metrizable apace $M$ is said to be a $Z$-set in $M$ \cite{vm:89},   if the set $C(Q,M\backslash F)$ is dense in $C(Q,M)$ with respect to the uniform convergence topology, where $Q$ denotes the Hilbert cube. 

\begin{lem}
Let $\alpha\colon X\to (0,\infty)$ be a positive continuous function, $V\subset\uin^{n+1}$ open and $g_0\in C(X,\uin^{n+1})$.
Then $\psi_V(y)\cap\overline{B}(g_0,\alpha)$ is a $Z$-set in $\overline{B}(g_0,\alpha)$ for every $y\in Y$, where $\overline{B}(g_0,\alpha)$ is considered as a subspace of $C(X,\uin^{n+1})$ with the uniform convergence topology. 
\end{lem}

\begin{proof}
The proof of this lemma follows very closely the proof of \cite[Lemma 2.8]{tv:01}. For sake of completeness we provide a sketch.
In this proof all function spaces are equipped with the uniform convergence topology generated by the Euclidean metric $d$ on $\uin^{n+1}$.
Since, by Lemma 2.4, $\psi_V$ has a closed graph, each $\psi_V(y)$
is closed in $\overline{B}(g_0,\alpha)$. We need to show that, for fixed $y\in Y$, $\delta>0$ and a map $u\colon Q \to \overline{B}(g_0,\alpha)$ there exists a map
$v\colon Q\to\overline{B}(g_0,\alpha)\backslash\psi_V(y)$ which is $\delta$-close to $u$. Observe first that $u$ generates $h\in C(Q\times X,\uin^{n+1})$, $h(z,x)=u(z)(x)$, such that
$d(h(z,x),g_0(x))\leq\alpha (x)$ for any $(z,x)\in Q\times X$. Since $f^{-1}(y)$ is compact, take $\lambda\in (0,1)$ such that $\lambda\sup\{\alpha (x):x\in f^{-1}(y)\}<\displaystyle\frac{\delta}{2}$ and define $h_1\in C(Q\times f^{-1}(y),\uin^{n+1})$ by $h_1(z,x)=(1-\lambda)h(z,x)+\lambda g_0(x)$. Then, for every $(z,x)\in Q\times f^{-1}(y)$, we have \\

\smallskip\noindent   
(2) \hbox{}~~~~~~$d(h_1(z,x),g_0(x))\leq (1-\lambda)\alpha (x)<\alpha (x)$ \\

\smallskip\noindent
and

\smallskip\noindent
(3) \hbox{}~~~~~~$d(h_1(z,x),h(z,x))\leq\lambda\alpha (x)<\displaystyle\frac{\delta}{2}$.

\smallskip\noindent
Let $\displaystyle q<\min\{r,\frac{\delta}{2}\}$, where $r=\inf\{\alpha (x)-d(h_1(z,x),g_0(x)):(z,x)\in Q\times f^{-1}(y)\}$. 
Since $\dim f^{-1}(y)\leq n$, by Lemma 2.3 (applied to the product $Q\times f^{-1}(y)$),   
there is a map $h_2\in C(Q\times f^{-1}(y),\uin^{n+1})$ such that $d(h_2(z,x),h_1(z,x))<q$ and $h_2(\{z\}\times f^{-1}(y))$ doesn't
contain $V$ for each $(z,x)\in Q\times f^{-1}(y)$. Then, by $(2)$ and $(3)$, for all $(z,x)\in Q\times f^{-1}(y)$ we have \\

\smallskip\noindent
(4) \hbox{}~~~~~~$d(h_2(z,x),h(z,x))<\delta$ and $d(h_2(z,x),g_0(x))<\alpha (x)$. \\

\smallskip\noindent
Because both $Q$ and $f^{-1}(y)$ are compact, $u_2(z)(x)=h_2(z,x)$ defines 
the map $u_2\colon Q\to C(f^{-1}(y),\uin^{n+1})$. 
Since the map $\pi\colon\overline{B}(g_0,\alpha)\to C(f^{-1}(y),\uin^{n+1})$, $\pi (g)=g|f^{-1}(y)$, is continuous and open (with respect to the uniform convergence topology), we can see that
$u_2(z)\in\pi (\overline{B}(g_0,\alpha))$ for every $z\in Q$ and
$\theta(z)=\overline{\pi^{-1}(u_2(z))\cap B_{\delta}(u(z))}$ defines a convex-valued map from $Q$ into $\overline{B}(g_0,\alpha)$ which is lower semi-continuous. By  the Michael selection theorem \cite[Theorem 3.2"]{em:56}, there is a continuous selection $v\colon Q\to C(X,\uin^{n+1})$ for $\theta$. Then $v$ maps $Q$ into $\overline{B}(g_0,\alpha)$ and $v$ is $\delta$-close to $u$. Moreover, for any $z\in Q$ we have $\pi(v(z))=u_2(z)$ and $V\not\subset u_2(z)( f^{-1}(y))$. Hence, 
$v(z)\not\in\psi_V(y)$ for any $z\in Q$, i.e.
$v\colon Q\to\overline{B}(g_0,\alpha)\backslash\psi_V(y)$.       
\end{proof}

We are now in a position to finish the proof of Theorem 1.1.

\begin{lem}
The set $\mathcal{H}_V$ is dense in $C(X,\uin^{n+1})$ for every open $V\subset\uin^{n+1}$.
\end{lem}

\begin{proof}
We need to show that, for fixed $g_0\in C(X,\uin^{n+1})$ and a continuous function $\alpha\colon X\to (0,\infty)$, there exists $g\in \overline{B}(g_0,\alpha)\cap\mathcal{H}_V$. The space $C(X,\uin^{n+1})$ with the uniform convergence topology 
is a closed convex subspace of the Banach space $E$ consisting of all bounded continuous maps from $X$ into $\R^{n+1}$.
We define the set-valued map $\phi$ from $Y$ into $C(X,\uin^{n+1})$, 
$\phi(y)=\overline{B}(g_0,\alpha)$, $y\in Y$. According to Lemma 2.5, 
$\overline{B}(g_0,\alpha)\cap\psi_V(y)$ is a $Z$-set in $\overline{B}(g_0,\alpha)$ for every $y\in Y$. So, we have a lower semi-continuous closed and convex-valued map $\phi\colon Y\to\mathcal{F}(E)$ and another map $\psi_V\colon Y\to\mathcal{F}(E)$ with a closed graph (see Lemma 2.4) such that 
$\phi(y)\cap\psi_V(y)$ is a $Z$-set in $\phi(y)$ for each $y\in Y$. Moreover, $Y$ is a $C$-space, so we can apply 
\cite[Theorem 1.1]{gv:99} to obtain a  continuous map $h\colon Y\to C(X,\uin^{n+1})$ with $h(y)\in\phi(y)\backslash\psi_V(y)$ for every $y\in Y$. Then
$g(x)=h(f(x))(x)$, $x\in X$, defines a map $g\in \overline{B}(g_0,\alpha)$. On the other hand, $h(y)\not\in\psi_V(y)$, $y\in Y$, implies that $g\in\mathcal{H}_V$.
 \end{proof} 

\noindent {\em Proof of Corollary $1.2$.}

As in the proof of Theorem 1.1, we can suppose that $f$ is perfect.
We first consider the case when $m$ is an integer $\geq n+1$.  Let $exp_{n+1}$ be the family of all subsets of $A=\{1,2,..,,m\}$ having cardinality $n+1$ and let $\pi_B\colon\uin^m\to\uin^B$ denote the corresponding projections, $B\in exp_{n+1}$.  It can be shown that $C(X,\uin^m)=C(X,\uin^B)\times C(X,\uin^{A\backslash B})$, so each  projection
$p_B\colon C(X,\uin^m)\to C(X,\uin^B)$ is open. Since, by Theorem 1.1, every set 
$\mathcal{H}_B=\{g\in C(X,\uin^B): \dim g(f^{-1}(y))\leq n\mbox{ for all $y\in Y$}\}$ is dense and $G_{\delta}$ in $C(X,\uin^B)$, so is the set $p_B^{-1}(\mathcal{H}_B)$ in $C(X,\uin^m)$. Consequently, the intersection $\mathcal{H}_m$ of all $\mathcal{H}_B$,
$B\in exp_{n+1}$, is also dense and $G_{\delta}$ in $C(X,\uin^m)$. Moreover, if $g\in\mathcal{H}_m$ and $y\in Y$, then  $\dim \pi_B(g(f^{-1}(y)))\leq n$ for any $B\in exp_{n+1}$.  The last inequalities, according to a result of N\"{o}beling \cite[Problem 1.8.C]{re:95}, imply $\dim g(f^{-1}(y)))\leq n$.

Now, let $m=\omega$ and $exp_{<\omega}$ denote the family of all finite sets $B\subset\omega$ of cardinality $|B|\geq n+1$.
Keeping the above notations, for any $B\in exp_{<\omega}$, $\pi_B\colon Q=\uin^\omega\to\uin^B$ and 
$p_B\colon C(X,Q)\to C(X,\uin^B)$ stand for the corresponding projections.  Then the intersection $\mathcal{H}_\omega$ of all 
$p_B^{-1}(\mathcal{H}_{B})$ is dense and $G_\delta$ in $C(X,Q)$.  We need only to check that $\dim g(f^{-1}(y))\leq n$ for any
$g\in\mathcal{H}_\omega$ and $y\in Y$. And this is certainly true, take an increasing sequence $\{B(k)\}$ in  $exp_{<\omega}$ which covers $\omega$ and consider the inverse sequence $\mathcal{S}=\{\pi_{B(k)}(g(f^{-1}(y))), \pi^{k+1}_k\}$, where
$\pi^{k+1}_k\colon \pi_{B(k+1)}(g(f^{-1}(y)))\to \pi_{B(k)}(g(f^{-1}(y)))$ are the natural projections.  Obviously, 
 $g(f^{-1}(y))$ is the limit space of $\mathcal{S}$. Moreover,  $g\in\mathcal{H}_\omega$ implies that 
$\pi_{B(k)}\circ g\in\mathcal{H}_{B(k)}$ for any $k$, so   all
 $\pi_{B(k)}(g(f^{-1}(y)))$ are at most $n$-dimensional. Hence,  $\dim g(f^{-1}(y))\leq n$.

\bigskip

\end{document}